\documentclass[12pt,a4paper]{article}
\usepackage[utf8]{inputenc}  
\usepackage{amssymb,amsmath,theorem}

 \let\kappa=\varkappa                  %
\newcommand\qed{\ifhmode\unskip\nobreak\fi\quad            %
   \ifmmode\square\else\hbox{$\square$}\fi}                %
\newcommand\pin{\kern.0833em}                              %

\theorempreskipamount=\medskipamount \theorempostskipamount=\medskipamount

 \newcommand\proof[1]{{\it Proof#1}\,}
 \newtheorem{theorem}{Theorem}
 \newtheorem{lemma}[theorem]{Lemma}
 \newtheorem{proposition}[theorem]{Proposition}
 
 {\theorembodyfont{\normalfont}
 \newtheorem{remark}{Remark}

\DeclareMathOperator{\Erg}{Erg}

\begin{document}

\begin{center}
 {\bf \large An ergodic support of a dynamical system and a natural\\[3pt]
 representation of Choquet distributions for invariant measures}

\bigskip\normalsize\rm

V.\,I.\ Bakhtin

\smallskip

{\it Belarusian State University, Belarus $($e-mail: bakhtin@tut.by\/$)$}
\end{center}

\vspace{-12pt}

\renewcommand\abstractname{}
\begin{abstract} \noindent
An ergodic support $X_0$ of a dynamical system $(X,T)$ with metrizable compact phase space $X$ is the set of all points $x\in X$ such that the corresponding sequence of empirical measures $\delta_{x,n} = (\delta_x +\delta_{Tx}+\dots +\delta_{T^{n-1}x})/n$ converges weakly to some ergodic measure. For every invariant probability measure $\mu$ on $X$ it is proven that $\mu(X_0) =1$ and Choquet distribution $\mu^*$ on the set of ergodic measures $\Erg X$ has the natural representation $\mu^*(A) =\mu(\{\pin x\in X_0 : \lim\delta_{x,n} \in A\pin\})$, where $A\subset \Erg X$.
\end{abstract}

\bigskip

{\textbf{Keywords:} {\itshape ergodic support, ergodic measure, Choquet distribution}

\medskip

}

\bigskip\medskip

Let $T\colon X\to X$ be a continuous mapping of a compact metric space $X$ into itself. Denote by $M(X)$ the set of all Borel probability measures on $X$, by $M_T(X)$ the subset of $T$-invariant measures from $M(X)$, and by $\Erg X$ all ergodic measures of the dynamical system $(X,T)$. Obviously,
\begin{equation*}
  \Erg X \subset M_T(X) \subset M(X).
\end{equation*}

\medskip

We suppose that $M(X)$ and $M_T(X)$ are equipped with a weak topology generated by continuous functions on $X$. Then by Alaoglu's theorem they are both compact. It is well known that the weak topology on $M(X)$ is metrizable via the metric given by formula
\begin{equation} \label{,,1}
 \rho(\nu,\mu) \pin=\, \sum_{i=1}^\infty \frac{1}{2^i}
 \bigg|\! \int_X f_i\,d\nu \,-\int_X f_i\,d\mu\pin \bigg|, \qquad \nu,\mu\in M(X),
\end{equation}
where the sequence of functions $f_i\in C(X)$ is dense in the unit ball in $C(X)$.

Below we will use the following statement proven in \cite[Proposition 1.3]{Felps}: \emph{the set\/ $\Erg X$ of all ergodic measures on a compact metric space\/ $X$ is a Borel set of type\/ $G_\delta$.}

Let $\delta_x$ be a unit measure supported at $x\in X$. Let us define a sequence of \emph{empirical measures} $\delta_{x,n}\in M(X)$ by the rule
\begin{equation} \label{,,2}
 \delta_{x,n} = \frac{\delta_{x}+ \delta_{Tx} + \dots+\delta_{T^{n-1}x}}{n}.
\end{equation}

For any fixed $n$ the empirical measure $\delta_{x,n}$ depends continuously on $x\in X$. In addition, for any fixed $x$ all accumulation points of the sequence $\delta_{x,n}$ are invariant measures, i.\,e., belong to $M_T(X)$ (see, for instance, \cite[Lemma 4]{B-S}).

It will be convenient for us to use the notation
\begin{equation*}
  \delta_{x,n}[f] = \frac{f(x)+f(Tx)+\dots+f(T^{n-1}x)}{n}, \qquad \mu[f] =\int_X f\,d\mu.
\end{equation*}

For each invariant measure $\mu\in M_T(X)$ its \emph{basin} is the set
\begin{equation*}
 B(\mu) = \{\pin x\in X :\pin \delta_{x,n} \to \mu\pin\},
\end{equation*}
where $\delta_{x,n} \to \mu$ means the weak convergence. Equivalently, the basin can be defined by means of the equality
\begin{equation}
  B(\mu) = \big\{ x\in X :\pin \delta_{x,n}[f_i]\to \mu[f_i],\ \ i\in\mathbb N \big\}, \label{,,3}
\end{equation}
where the functions $f_i\in C(X)$ are taken from \eqref{,,1}.

\pagebreak[2]

\begin{proposition} \label{..1}
The set\/ $B(\mu)$ is\/ $T$-invariant and a Borel one.
\end{proposition}

\proof. From \eqref{,,2} it is seen that the sequences $\delta_{x,n}$ and $\delta_{Tx,n}$ have the same limits. This implies the invariance of $B(\mu)$. The set $B(\mu)$ can be defined via operations of countable union and countable intersection of Borel sets by the formula
\begin{equation*}
  B(\mu) \pin =\pin
  \bigcap_{k=1}^\infty \bigcup_{N=1}^\infty \bigcap_{n=N}^\infty
  \left\{x\in X :\pin \rho(\delta_{x,n},\mu) <\frac{1}{k} \right\}.
\end{equation*}
Hence it is Borel. \qed

\medskip

An \emph{ergodic support} of the dynamical system $(X,T)$ we call the set
\begin{equation*}
 X_0 = \bigcup\, \{\pin B(\mu) :\pin \mu\in \Erg X\pin\}.
\end{equation*}

In addition, let
\begin{gather*}
 X_1 = \bigcup\, \{\pin B(\mu) :\pin \mu\in M_T(X) \setminus \Erg X\pin\}, \\[3pt]
 X_2 = \{\pin x\in X :\pin \nexists\, \lim \delta_{x,n}\pin\} = X \setminus (X_0\cup X_1).
\end{gather*}

\begin{proposition} \label{..2}
The sets\/ $X_0,\, X_1,\, X_2$ are\/ $T$-invariant and Borel.
\end{proposition}

\proof. The invariance of $X_0,\, X_1,\, X_2$ follows from the invariance of basins $B(\mu)$.

The set $X_2$ is Borel since it can be defined by the formula
\begin{equation*}
 X_2 \pin=\pin \bigcup_{k=1}^\infty \bigcap_{N=1}^\infty \bigcup_{n,m\ge N}
 \left\{ x\in X :\pin \rho(\delta_{x,n},\delta_{x,m}) \ge \frac{1}{k} \right\}.
\end{equation*}
Therefore the union $X_0\cup X_1$ is Borel as well.

Since the set $\Erg X$ is $G_\delta$ its complement $M_T(X)\setminus \Erg X$ is $F_\sigma$, that is, a countable union of closed sets:
\begin{equation*}
 M_T(X)\setminus \Erg X \pin=\, \bigcup_{i=1}^\infty F_i, \qquad F_i \ \, \text{is closed}.
\end{equation*}
Therefore $X_1$ may be presented as a countable union
\begin{equation*}
 X_1 \pin=\pin \bigcup_{i=1}^\infty X(F_i), \quad \text{where}\ \
 X(F_i) \pin=\, \bigcup\, \{\pin B(\mu) :\pin \mu\in F_i\pin\}.
\end{equation*}
Each $X(F_i)$ is Borel because it can be defined by the formula
\begin{equation*}
 X(F_i) \pin=\pin \bigcap_{k=1}^\infty \bigcup_{N=1}^\infty \bigcap_{n=N}^\infty
 \left\{ x\in X_0\cup X_1 :\pin \rho(\delta_{x,n},F_i) <\frac{1}{k}\right\}.
\end{equation*}
It follows that the sets $X_1$ and $X_0 =X\setminus (X_1\cup X_2)$ are Borel as well. \qed

\medskip

Next we present the following modification of the Birkhoff ergodic theorem.

\begin{theorem} \label{..3}
For each ergodic measure\/ $\mu\in \Erg X$ we have\/ $\mu(B(\mu)) =1$.
\end{theorem}

\proof. According to \eqref{,,3} the basin $B(\mu)$ can be defined as a countable intersection of the sets $\{\pin x\in X:\pin \delta_{x,n}[f_i] \to \mu[f_i]\pin\}$, where the collection of functions $f_i$ is dense in $C(X)$. By the ergodic theorem each of these sets has full measure $\mu$. Then its intersection has also full measure $\mu$. \qed

\medbreak

Notice that in the case of nonmetrizable compact space $X$ Theorem \ref{..3} is not valid.

A characteristic feature of the ergodic support $X_0$ of the dynamical system $(X,T)$ is that it bears all invariant measures. It is formulated in the next theorem.

\begin{theorem} \label{..4}
The equality\/ $\mu(X_0) =1$ holds for all\/ $\mu\in M_T(X)$.
\end{theorem}

Recall that $X$ is supposed to be a compact metric space here.

\medskip

\proof. Let us exploit the notation
\begin{equation*}
  f(\mu) =\int_X f\,d\mu, \qquad \mu\in M(X).
\end{equation*}
By Choquet's theorem \cite[page 14]{Felps} for any $\mu\in M_T(X)$ there exists a Borel probability measure $\mu^*$ on the set $\Erg X$ satisfying the identity
\begin{equation}\label{,,4}
  f(\mu) = \int_{\Erg X} f(\nu)\, d\mu^*(\nu), \qquad f\in C(X)
\end{equation}
(wherein $\mu^*$ is called the \emph{Choquet distribution} for $\mu$ and $\mu$ is called a \emph{barycenter} of $\mu^*$). By means of limiting process this identity extends to bounded Borel measurable functions (see Lemma~\ref{..5} below). Applying \eqref{,,4} to the characteristic function $f =I_{X_0}$ and exploiting Theorem \ref{..3} we obtain
\begin{equation*}
  \mu(X_0) =I_{X_0}(\mu) =\int_{\Erg X} I_{X_0}(\nu)\, d\mu^* =\int_{\Erg X} \nu(X_0)\, d\mu^* =1. \qed
\end{equation*}

\begin{lemma} \label{..5}
Equality\/ \eqref{,,4} holds true for all bounded Borel functions\/ $f$.
\end{lemma}

\proof. Let us use a modified version of the construction procedure for classes of Baire functions from \cite[\S\medspace 43]{Hausdorff}, wherein we only replace monotone functional limits by uniformly bounded functional limits (in order to guarantee a passage to the limit in \eqref{,,4}).

Let $\omega_1$ be a minimal uncountable ordinal number. Set $\mathcal B_1 =C(X)$. By induction for every ordinal $\alpha <\omega_1$ define $\mathcal B_\alpha$ as the class of all bounded functions on $X$ that are pointwise limits of uniformly bounded sequences of functions from $\bigcup \{\pin\mathcal B_\beta : \beta <\alpha\pin\}$. Standard inductive reasoning shows that each class $\mathcal B_\alpha$ consists of Borel functions and forms an algebra. The union $\mathcal B =\bigcup\pin \{\pin \mathcal B_\alpha :\pin \alpha<\omega_1\pin\}$ possesses the same properties.

Now we check that if a uniformly bounded sequence of functions $f_i\in \mathcal B$ converges to a function $f$ then $f\in \mathcal B$. By construction $f_i\in \mathcal B_{\alpha_i}$ for some $\alpha_i <\omega_1$. Set $\alpha =\sup_i \alpha_i$. This $\alpha$ is a countable ordinal number and hence $f\in \mathcal B_{\alpha+1}\subset \mathcal B$.

It is easily seen that $\mathcal B$ contains characteristic functions of open and closed subsets of $X$. On the other hand, the subsets of $X$ whose characteristic functions belong to $\mathcal B$ compose a $\sigma$-algebra. This $\sigma$-algebra contains all Borel sets, and hence $\mathcal B$ contains all bounded Borel functions (and nothing else).

Applying the transfinite induction over $\alpha<\omega_1$ and Lebesgue's dominated convergence theorem, we prove validity of \eqref{,,4} for all functions $f\in \mathcal B_\alpha$ and, as a result, for all bounded Borel functions. \qed

\begin{remark}
As an extra bonus, the same reasoning shows that for any bounded Borel function $f$ on $X$ the corresponding function $f(\nu)$ on $\Erg X$ is also Borel.
\end{remark}

\begin{remark}
In the case of nonmetrizable topological space $X$ the algebra $\mathcal B$ built in the proof of Lemma \ref{..5} consists of all bounded Baire functions. In the case of metrizable $X$ the notions of Borel functions and Baire functions coincide.
\end{remark}

Consider a mapping $V\colon X_0\to \Erg X$ defined by the formula
\begin{equation} \label{,,5}
  V(x) =\lim_{n\to\infty} \delta_{x,n}, \qquad x\in X_0
\end{equation}
(of course, for $x\in X_1$ this formula also makes sense and gives $V(x)\in M_T(X)$, but this is not interesting for us).

\begin{proposition} \label{..6}
The mapping defined in\/ \eqref{,,5} is Borel measurable.
\end{proposition}

\proof. For every $f\in C(X)$ restriction of the function $f(V(x))$ to $X_0$ is Borel measurable because it can be represented as a limit of a sequence of continuous functions:
\begin{equation*}
  f(V(x)) = \lim_{n\to \infty} f(\delta_{x,n}) =
  \lim_{n\to\infty} \frac{f(x)+f(Tx)+\dots+f(T^{n-1}x)}{n}.
\end{equation*}
It follows that the mapping \eqref{,,5} is Borel measurable. \qed

\medskip

By virtue of Proposition \ref{..6} the mapping \eqref{,,5} transfers Borel measures $\mu$ supported on $X_0$ to Borel measures $\mu^* =V(\mu)$ on the set $\Erg X$ by the rule
\begin{equation} \label{,,6}
 \mu^*(A) =\mu(V^{-1}(A)), \qquad A\subset \Erg X.
\end{equation}
It turns out that for any measure $\mu\in M_T(X)$ the corresponding measure $\mu^* =V(\mu)$ on $\Erg X$, defined in \eqref{,,6}, is Choquet distribution from \eqref{,,4}. This follows from the next theorem.

\begin{theorem} \label{..7}
For each\/ $\mu\in M_T(X)$ the measure\/ $\mu^* =V(\mu)$ is Choquet distribution on\/ $\Erg X$ with barycenter\/ $\mu$, and for all bounded Borel functions\/ $f$ on\/ $X$ the following equalities hold true\/$:$
\begin{equation} \label{,,7}
 f(\mu) =\int_{X_0} f(V(x))\, d\mu(x) =\int_{\Erg X} f(\nu)\, d\mu^*(\nu).
\end{equation}
\end{theorem}

\proof. Since $\mu(X_0) =1$ the measure $\mu^* =V(\mu)$ is probability. For $f\in C(X)$ the left-hand equality in \eqref{,,7} comes out in the limit from the equality
\begin{equation*}
  \int_{X_0} f(x)\,d\mu(x) =\int_{X_0} \frac{f(x) +f(Tx) +\dots +f(T^{n-1}x)}{n}\, d\mu(x),
\end{equation*}
and the right-hand equality in \eqref{,,7} follows from the definition of $\mu^*$. After that \eqref{,,7} extends to all bounded Borel functions in the same way as in the proof of Lemma \ref{..5}. \qed

\begin{theorem} \label{..8}
If\/ $\tau(\mu)$ is an affine upper semicontinuous function on\/ $M_T(X)$ then
\begin{equation} \label{,,8}
  \tau(\mu) \pin= \int_{\Erg X} \tau(\nu)\,d\mu^*(\nu) \pin=\int_{X_0} \tau(V(x))\,d\mu(x),
  \qquad \mu\in M_T(X).
\end{equation}
\end{theorem}

\proof. Since $\mu$ is a barycenter of the measure $\mu^* =V(\mu)$ the left-hand equality in \eqref{,,8} follows from \cite[Lemma 10.7]{Felps}. The right-hand equality in \eqref{,,8} follows from the definition of the measure $\mu^*$. \qed

\begin{remark}
The author knows an example of a bounded affine Borel function $\tau(\mu)$ on the set $M_T(X)$ that breaks the left-hand equality in \eqref{,,8}.
\end{remark}



\end{document}